# Luenberger-Type Cubic Observers for State Estimation of Linear Systems

Mohammad Mahdi Share Pasand[1]

[1] Department of Electrical Engineering, Faculty of Technology and Engineering, Standard Research Institute, Alborz, Iran, PO Box 31585-163, Email: sharepasand@standard.ac.ir

*Abstract.* This paper introduces a new nonlinear observer for state estimation of linear time invariant systems. The proposed observer contains a (nonlinear) cubic term in its error dynamics. The cubic observer can yield improved performance while possessing desired robustness properties. Unlike previously proposed observers, estimation error dynamics of the cubic observer are nonlinear. Convergence criteria, robustness properties and observer based feedback control with cubic observers are addressed. Simulation examples are included as well which show significant performance improvement compared to linear observers.

## 1. Introduction

Since their introduction, observers have been extensively studied and used in the literature for both linear and nonlinear dynamical systems [1-3]. In practice, not all state variables can be directly measured by sensors. As a result, in many applications requiring the whole state variable vector, it is inevitable to use an observer to construct the state variable from the measured output [1]. Observers may be used for observer-based feedback control [4-8], fault detection [9-14] and parameter estimation [1]. Unknown or partially known input [3, 10], robust [7, 15], optimal [2, 15] and functional [16] observers have been investigated in the literature for linear [16] as well as descriptor [3, 16], delayed [9], time varying [2], linear parameter varying [13, 14] and nonlinear systems [1, 4-8, 11, 12]. The literature of robust and optimal linear observers for linear systems is very mature [15]. Nonlinear state estimation including sliding mode and high gain observers have received research attention as well [7, 11]. A significantly different approach to observer design has been mentioned in [17] where the state trajectory of the so-called interval observer, does not need to converge to the plant state, but is desired to be in a pre-defined vicinity of it. Not only linear ordinary observers, but also delayed [18], switching [19], descriptor [10] and nonlinear [1] observers have been practiced in the literature. Finite time [21], perfect [22] and proportional, integral, derivative [23-25] observers have been introduced as well to enhance performance of linear observers in different manners.

Though the observer model has become more and more complicated to imitate complicated plant dynamics, it is common to keep the estimation error dynamics linear (or delayed linear [16]) to facilitate stability analysis. As a result, less attention has been given to nonlinear estimation error dynamics. In addition, to the best of the author's knowledge, use of nonlinear observers for linear systems has not been considered yet. This paper suggests a nonlinear observer, namely the cubic observer, as a generalization to linear observers, in order to improve the observer performance while maintaining desired robustness properties. The proposed cubic observer structure resembles the structure of Luenberger observers. The motivation of introducing the cubic observer is to achieve a fast response with smaller error norm. Achieving fast observer responses has been the aim of several researches [26-33] and is essential to observer based feedback control [4], finite time [21] and time-optimal [34] observers as well as fault detection observers [2, 3, 10]. Compared to a similar linear observer, the cubic observer yields a faster response in the beginning of the observation course and a faster descending Lyapunov function when error norms are equal. The cubic term can be made arbitrarily small or set to zero, to transform the cubic observer to a linear one. Therefore, a cubic observer is rather a generalization than an alternative to the linear observers. The cubic observer requires generic existence conditions and has a simple design method. Unknown and partially known input [3], proportional-integral-derivative [23-25] and optimal cubic observers are among future possible works.

This paper is organized in four sections. In the second section, main results are given. In section 2.1., the cubic observer is introduced and its global stability conditions are discussed. In addition, we show that the cubic observer has performance advantages over a similar linear observer. In section 2.2., it is shown that cubic observers possess some robustness properties. Section 2.3., provides conditions under which the cubic observer can be deployed in combination with linear state feedback for regulation and control. It is shown that compared to linear observers, no additional requirement is needed for observer based feedback control with cubic observers. The third section provides examples. The fourth section concludes the paper.

## 2. Luenberger-Type Cubic Observers

### 2.1. Convergence Criteria and Performance Advantages

Assume a linear time invariant system described by state space equations (1).

$$\dot{x}(t) = Ax(t) + Bu(t), \quad y(t) = Cx(t) \tag{1}$$

Vector $x(t) \in R^n$ represents the state variable, $y(t) \in R^{n_y}$ is the measured output and matrices $A_{n \times n}$, $B_{n \times n_u}$ and $C_{n_y \times n}$ are the state, input and output matrices respectively. Vector $u(t) \in R^{n_u}$ is a known (measured) exogenous input. It is assumed throughout the paper that the pair $(A, C)$ is observable. In order to estimate the state of system (1) from the measured output, it is common to use a linear observer (2).

$$\dot{\hat{x}}_l(t) = (A - LC)\hat{x}_l(t) + Bu(t) + Ly(t) \tag{2}$$

In (2), $\hat{x}_l(t) \in R^n$ represents the estimated state. Matrix $L_{n \times n_y}$ is the observer gain. Subscript $l$ is used to emphasis that observer (2) is a linear dynamical system. The estimation error is defined as:

$$e_l(t) = x(t) - \hat{x}_l(t) \tag{3}$$

This error is governed by the following dynamical equations:

$$\dot{e}_l(t) = \dot{x}_n(t) - \dot{\hat{x}}_l(t) = (A - LC)e_l(t) \tag{4}$$

Now assume a cubic observer of the form (5).

$$\dot{\hat{x}}_c(t) = (A - L_c C)\hat{x}_c(t) + L_c y(t) + Bu(t) - e_c^T(t) C^T \theta C e_c(t) N_c C e_c(t) \tag{5}$$

In which $\hat{x}_c(t) \in R^n$ is the estimated state vector for the cubic observer which is aimed to converge to the state vector of (1). Matrices $L_c \in R^{n \times n_y}$, $N_c \in R^{n \times n_y}$ and $\theta \in R^{n_y \times n_y}$ are the observer parameter matrices and the estimation error is defined as follows:

$$e_c(t) = x(t) - \hat{x}_c(t) \tag{6}$$

The estimation error dynamics are:

$$\dot{e}_c(t) = (A - L_c C)e_c(t) + e_c^T(t) C^T \theta C e_c(t) N_c C e_c(t) \tag{7}$$

Observer structure (5) differs from most of previously proposed observers [1, 6-8] in the sense that it leads to nonlinear estimation error dynamics (7). When the estimation error shrinks as $t \to \infty$, the cubic term in (7) vanishes faster than the linear term. Therefore, the cubic observer (5) has a behavior close to a linear observer (2) in the vicinity of the origin $e_c(0) = e_l(0) = 0$. As evident from (4) and (7), the presence of the exogenous input does not affect the estimation error dynamics of linear or cubic observers. Therefore, dropping the exogenous input from (1) does not affect the results except for that of section 2.3.

**Remark 1.** Matrix $\theta$ determines the relative importance of each output. If $\theta = diag(r_1, \ldots r_{n_y})$, each $r_i$ will represent an importance index for the $i^{th}$ output. Another approach is to determine $\theta$ based on the role of each output on the estimation cost after a linear observer is deployed for state estimation. For this purpose, we may first deploy a linear observer and then determine $\theta$ based on the observer response.

Defining the Lyapunov function candidate $V_c = e_c^T P e_c$, in which $P \in R^{n \times n}$ is to be determined, a stability criterion is derived as (8).

$$e_c^T\left((A - L_c C)^T P + P(A - L_c C)\right)e_c + (e_c^T C^T \theta C e_c) \times e_c^T (P N_c C + C^T N_c^T P) e_c < 0 \qquad (8)$$

We choose $L_c$ such that $A - L_c C$ is Hurwitz, satisfying the Lyapunov equation (9).

$$(A - L_c C)^T P + P(A - L_c C) = -Q < 0 \qquad (9)$$

It can be shown that for any $Q = Q^T > 0$, there exists a unique solution $P = P^T > 0$ to (9) if the pair $(A, C)$ is observable [35]. If (9) holds, (8) will be fulfilled if:

$$(e_c^T C^T \theta C e_c) \times e_c^T(t)(P N_c C + C^T N_c^T P) e_c(t) < 0 \qquad (10)$$

If $\theta \geq 0$ and $P N_c C + C^T N_c^T P < 0$, (10) will be fulfilled. *Theorem 1* summarizes the above discussion.

**Theorem 1.** If the positive definite symmetric matrix $P \in R^{n \times n}$ and $C^T \theta C \geq 0$ satisfy (11)-(13), the estimation error dynamics of observer (5) for system (1) will be globally stable.

$$\nexists\, v \neq 0; (A - L_c C)v + v^T C^T \theta C v N_c C v = 0 \qquad (11)$$

$$\exists\, P = P^T > 0, N_c, L_c \,; \begin{cases} (A - L_c C)^T P + P(A - L_c C) < 0 & (12) \\ P N_c C + C^T N_c^T P < 0 & (13) \end{cases}$$

Matrices $\theta = \theta^T \geq 0$, $N_c$ and $L_c$ are to be chosen to fulfill (11)-(13).

**Proof:** condition (11) guarantees that the origin is the only equilibrium for (7), while (12)-(13) guarantee that the Lyapunov function derivative given in (8)-(10) is negative. The proof is complete. ∎

**Lemma 1.** Assume a non-symmetric, real valued, square matrix $M \in R^{q \times q}$. If the matrix $M + M^T$ is negative definite, then $M$ is also negative definite in the sense that it fulfills:

$$v^T M v < 0 \,;\, \forall v \in R^q$$

**Proof:** Due to assumption $M + M^T < 0$, we have:

$$v^T M v + v^T M^T v < 0 \,;\, \forall v \in R^q$$

Since the two terms in the left hand side of the inequality are scalars and transposes of each other, they are equal. Therefore;

$$v^T M v + v^T M^T v = 2v^T M v < 0 \,;\, \forall v \in R^q$$

This completes the proof. ∎

**Remark 2.** *Theorem 1* can be alternatively established with Lyapunov function $V_{cz} = 1 - \exp(-e_c^T P e_c)$ and the Zubov's equation [34], for which a convergence region is derived as: $0 < V_{cz} < 1$.

**Theorem 2.** For a given $C^T \theta C \geq 0$, a solution to (11) and (13) is given by (14). (Scalar $\gamma > 0$ is arbitrary.)

$$N_c = -\gamma P^{-1} C^T \theta \tag{14}$$

*Proof:* The proof follows by defining the Lyapunov candidate function $V = e_c^T P e_c$. Due to observability assumption for the pair $(A, C)$ one can choose $L_c$ to make $(A - L_c C)$ stable. By assigning (14), condition (11) will be fulfilled. To show this, assume that there exists a non-zero equilibrium for (7). i.e.:

$$\exists \, v^* \neq 0; (A - L_c C) v^* + v^{*T} C^T \theta C v^* N_c C v^* = 0 \tag{15}$$

The scalar $v^{*T} C^T \theta C v^*$ is non-zero. To show this, notice that if this scalar is zero, then (15) requires that $(A - L_c C) v^* = 0$ with a nonzero $v^*$. This is not possible since $(A - L_c C)$ is assumed to be stable and thus of full rank, having an inverse matrix. Multiplying both sides of (15) by $(A - L_c C)^{-1}$, (16) will be obtained.

$$v^* = -v^{*T} C^T \theta C v^* (A - L_c C)^{-1} N_c C v^* \tag{16}$$

Left multiply the both parts of this equation by the vector $v^{*T} C^T \theta C$ to obtain:

$$v^{*T} C^T \theta C v^* = -v^{*T} C^T \theta C v^* v^{*T} C^T \theta C \, (A - L_c C)^{-1} N_c C v^* \tag{17}$$

Dividing by the non-zero scalar $v^{*T} C^T \theta C v^*$, it is resulted that:

$$-1 = v^{*T} C^T \theta C \, (A - L_c C)^{-1} N_c C v^*$$

This equation can't be fulfilled if the matrix $C^T \theta \, C (A - L_c C)^{-1} N_c C$ is positive semi-definite i.e. if:

$$v^T C^T \theta C \, (A - L_c C)^{-1} N_c C v \geq 0 \,;\, \forall v \in R^n$$

Substituting (14):

$$C^T \theta C \, (A - L_c C)^{-1} N_c C = -\gamma C^T \theta C \, (A - L_c C)^{-1} P^{-1} C^T \theta \, C$$

Note that:

$$P(A - L_c C) \leq 0 \Leftrightarrow -\gamma C^T \theta C \, (A - L_c C)^{-1} P^{-1} C^T \theta \, C \geq 0 \tag{18}$$

Since $L_c$ is chosen to stabilize $(A - L_c C)$, it is guaranteed that: $(A - L_c C)^T P + P(A - L_c C) < 0$. Therefore, according to *Lemma 1*, matrix $P(A - L_c C)$ is negative definite. This completes the proof. ∎

Scalar $\gamma > 0$ can be used for optimization while stability of observer is assured via *Theorem 2*. In the remainder of this section, we provide comparisons between the cubic and linear observers.

***Theorem 3.*** Denote the Lyapunov functions of the linear and cubic estimation error dynamics (4) and (7) by $V_l(t)$ and $V_c(t)$ respectively. If (11)-(13) are fulfilled for observer (5), $Q = \rho I, \rho > 0$ and $L = L_c$, then (19) holds.

$$e_c^T(t)e_c(t) = e_l^T(t)e_l(t) \Rightarrow \dot{V}_c(t) \leq \dot{V}_l(t) \tag{19}$$

***Proof.*** Write the Lyapunov derivative as (8). If (13) holds;

$$\dot{V}_c(e_c) \leq -e_c^T(t)Qe_c(t) = -\rho e_c^T(t)e_c(t)$$

Also note that;

$$\dot{V}_l(t) = -e_l^T(t)Qe_l(t) = -\rho e_l^T(t)e_l(t)$$

The proof completes here. ∎

*Theorem 3* alongside with equation (8), states an advantage of the cubic observer. When the estimation error norm is larger for the cubic observer, its Lyapunov function will descend faster than a linear observer as can be seen from (8). Furthermore, when the two observers have estimation errors with equal norms, the cubic observer has a faster response. The following *Remark* states another advantage of the cubic observer.

***Remark 3.*** If (11)-(13) are fulfilled for observer (5), $L = L_c$ and $e_c(0) = e_l(0)$, then based on (10), we have $\dot{V}_c(0) \leq \dot{V}_l(0)$. The larger the initial estimation error $e_c(0) = e_l(0)$ is, the larger-negative is the cubic term (10). Therefore, the estimated state of cubic observer takes a much larger step towards the actual state, if the initial discrepancy is larger.

## 2.2. Robustness Properties

In this section, robustness properties of cubic observers are addressed. Consider the uncertain system (20).

$$\dot{x}(t) = A(\varepsilon)x(t) + Bu(t), \qquad y(t) = Cx(t) \tag{20}$$

Assume that $(A(\varepsilon), C)$ fulfills the quadratic stabilizability property (21).

$$\exists P = P^T > 0 \ \& \ L_c; \ \ P(A(\varepsilon) - L_c C) + (A(\varepsilon) - L_c C)^T P < 0 \ ; \ \forall \varepsilon \tag{21}$$

***Theorem 4.*** Assume $L_c = L$ fulfills (11) and (21) for all uncertainties captured by $A(\varepsilon)$. The cubic observer (5) yields robustly stable estimation error dynamics for all uncertainties captured by $A(\varepsilon)$ if (13) holds.

***Proof:*** According to assumption (11), the estimation error dynamics possess only one equilibrium point at the origin. Define the Lyapunov function candidate $V_c = e_c^T P e_c$. Take derivatives to obtain:

$$\dot{V}_c = \dot{V}_1 + \dot{V}_2$$

$$\dot{V}_1 = e_c^T \big( (A(\varepsilon) - L_c C)^T P + P(A(\varepsilon) - L_c C) \big) e_c , \qquad \dot{V}_2 = e_c^T C^T \theta C e_c e_c^T (PN_c C + C^T N_c^T P) e_c$$

The term $\dot{V}_1$ is negative due to (21). The term $\dot{V}_2$ is negative due to (13). The proof is complete. ∎

Note that $\dot{V}_2$ introduced in the proof of *Theorem 4,* is not affected by uncertainty and therefore provides a confidence margin for the cubic observer convergence when the estimation error of (20) is large.

Condition (21) may not be always satisfied. The following theorem states an alternative robustness property. This theorem also holds for linear observers as can be seen from the proof.

**Theorem 5.** Assume the linear uncertain system described by (20) fulfills (22).

$$A(\varepsilon) = A + \varepsilon I \quad , \quad \varepsilon_{min} \leq \varepsilon \leq \varepsilon_{max} \quad (22)$$

In which $\varepsilon$ is a scalar real valued perturbation with $\varepsilon_{max} \geq 0$ and $\varepsilon_{min} \leq 0$.[1] The estimation error dynamics of the cubic observer are robustly stable if (23) holds. (Matrices $P, Q$ are given in (9)).

$$\varepsilon_{max} \leq \frac{\lambda_{min}(Q)}{2\lambda_{max}(P)} \quad (23)$$

*Proof:* The proof follows by finding upper limits for the cubic observer estimation error dynamics with the perturbation $\varepsilon$. The Lyapunov function derivative fulfills:

$$\dot{V}_c = \dot{V}_{cN} + \dot{V}_{c\varepsilon}$$

$$\dot{V}_{cN} = e_c^T((A - L_cC)^T P + P(A - L_cC) + e_c^T C^T \theta C e_c (PN_cC + C^T N_c^T P))e_c \quad , \quad \dot{V}_{c\varepsilon} = 2\varepsilon e_c^T P e_c$$

The first term in the right hand side is upper-bounded by:

$$\dot{V}_{cN} \leq -\lambda_{min}(Q) e_c^T e_c$$

In which $\lambda_{min}(.)$ represents the minimum eigenvalue of a matrix. Note that matrix $Q$ is assumed symmetric positive definite. Therefore, all of its eigenvalues are real-valued positive. The second term fulfills:

$$\dot{V}_{c\varepsilon} \leq 2\varepsilon_{max}\lambda_{max}(P) e_c^T e_c$$

In which $\lambda_{max}(.)$ denotes the maximum eigenvalue of a matrix. Again, due to the fact that $P$ is positive definite symmetric, all of its eigenvalues are real-valued positive. Using the above inequality we obtain:

$$\dot{V}_c \leq (-\lambda_{min}(Q) + 2\varepsilon_{max}\lambda_{max}(P)) e_c^T e_c$$

In order for the Lyapunov derivative to be negative, it is sufficient that:

$$-\lambda_{min}(Q) + 2\varepsilon_{max}\lambda_{max}(P) < 0$$

This completes the proof. ∎

### 2.3. Observer-based Feedback Control

To proceed to the next theorem, assume the estimated state is fed to the input via a state feedback gain $K \in R^{n_u \times n}$ as in (24).

$$u(t) = -K\hat{x}_c(t) \quad (24)$$

---

[1] Assuming $\varepsilon_{max} \geq 0$ and $\varepsilon_{min} \leq 0$ does not compromise generality of the result. If both bounds are positive, we can replace $\varepsilon$ by $\varepsilon' = \varepsilon - \varepsilon_{min}$ with $\varepsilon'_{min} = 0$, $\varepsilon'_{max} = \varepsilon_{max} - \varepsilon_{min}$ and $A$ by $A + \varepsilon_{min}I$. Similar arguments hold if both bounds are negative.

**Theorem 6.** Assume that matrices $P_{n \times n} = P^T > 0$ and $C^T \theta C \geq 0$ satisfy (11)-(13). The closed loop control system governed by (1), (5) and (24) is stable if the positive definite symmetric matrices $P, P_1 \in R^{n \times n}$ fulfill (25).

$$\begin{bmatrix} (A - BK)^T P_1 + P_1(A - BK) & P_1 BK \\ K^T B^T P_1 & (A - L_c C)^T P + P(A - L_c C) \end{bmatrix} < 0 \qquad (25)$$

*Proof:* Form the closed loop system dynamics by combining (1), (5) and (24).

$$\begin{cases} \dot{x}(t) = (A - BK)x(t) + BKe_c(t) \\ \dot{e}_c(t) = (A - L_c C)e_c(t) + e_c^T(t)C^T \theta C e_c(t) N_c C e_c(t) \end{cases}$$

Defining the Lyapunov candidate function $V(e_c, x) = x^T P_1 x + e_c^T P e_c$ and taking derivatives we obtain:

$$\dot{V}(e_c, x) = x^T \big( (A - BK)^T P_1 + P_1(A - BK) \big) x + e_c^T \big( (A - L_c C)^T P + P(A - L_c C) \big) e_c + \cdots$$

$$\cdots e_c^T C^T \theta C e_c e_c^T (P N_c C + C^T N_c^T P) e_c + e_c^T K^T B^T P_1 x + x^T P_1 BK e_c = \cdots$$

$$[x^T \quad e_c^T] \Psi \begin{bmatrix} x \\ e_c \end{bmatrix} + e_c^T C^T \theta C e_c e_c^T (P N_c C + C^T N_c^T P) e_c$$

In which $\Psi$ is given in (25). Note that stability of $(A - BK)$ and $(A - L_c C)$ is necessary for (25). ∎

**Corollary 1.** The closed loop control system described by (1), (5) and (24) (i.e. linear feedback with a cubic observer) is stable if the closed loop control system described by (1), (2) and (24) (i.e. linear feedback with linear observer) is stable.

*Proof:* Rewrite left hand side of (25) as:

$$\begin{bmatrix} (A - BK)^T P_1 + P_1(A - BK) & P_1 BK \\ K^T B^T P_1 & (A - LC)^T P + P(A - LC) \end{bmatrix} = A_a^T P_a + P_a A_a$$

In which:

$$A_a = \begin{bmatrix} A - BK & BK \\ 0 & A - LC \end{bmatrix}, \qquad P_a = \begin{bmatrix} P_1 & 0 \\ 0 & P \end{bmatrix}$$

Dynamics of the closed loop system with linear feedback and linear observer can be written as $[\dot{x}^T \ \dot{e}_l^T]^T = A_a [x^T \ e_l^T]^T$. Note that $P_a > 0$ if and only if $P_1, P > 0$. The left side of (25) is the Lyapunov function derivative for (1) with observer (2) and state feedback (24). The proof is complete. ∎

## 3. Simulation examples

*Example 1.* In this example, we study the cubic observer performance for an integrator system with the following state space matrices.

$$A = \begin{bmatrix} 0 & 1 \\ 0 & 0 \end{bmatrix}, \quad B = \begin{bmatrix} 0 \\ 1 \end{bmatrix}, \quad C = [1 \ 0], \quad u(t) = \sin(t), \quad x_1(0) = x_2(0) = -3$$

The plant is observable. It is desired to estimate the second state from the measured output. Locating the linear observer poles at $[-2 \ -5]$ and using the pole placement technique, the linear observer gain is obtained as $L = [7 \ 10]^T$. Assuming $Q = 10 I_2$, the solution to the Lyapunov equation (12) is derived as:

$$P = \begin{bmatrix} 7.8571 & -5 \\ -5 & 4.2857 \end{bmatrix}$$

The cubic observer parameters are as follows:

$$L_c = L, \quad \theta = 10, \quad \gamma = 2, \quad N_c = -\gamma P^{-1} C^T \theta = \begin{bmatrix} -9.8824 \\ -11.5294 \end{bmatrix}$$

**Fig.1** shows the estimated and actual values for the second state with the linear and cubic observers in the first four seconds of the simulation course. The cubic observer provides a faster approach towards the actual state. **Fig.2** depicts the Lyapunov functions $V_c(t), V_l(t)$ in the first instants of the estimation for the two observer estimation errors introduced in (19). The faster decent of the cubic observer Lyapunov function is evident. **Fig.3** shows the two Lyapunov functions is a longer time window. It can be seen that when the estimation error is larger (due to initial conditions), the Lyapunov function for the cubic observer estimation error descends to large negative much faster. The decent rate approaches to zero eventually.

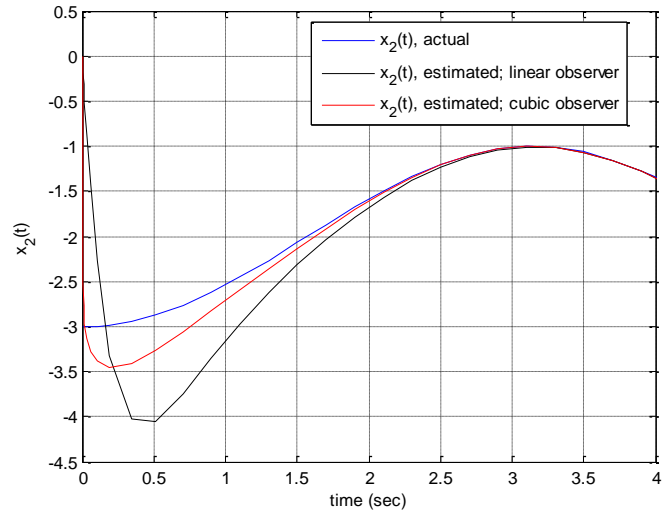

**Figure 1.** Estimated state with cubic and linear observers

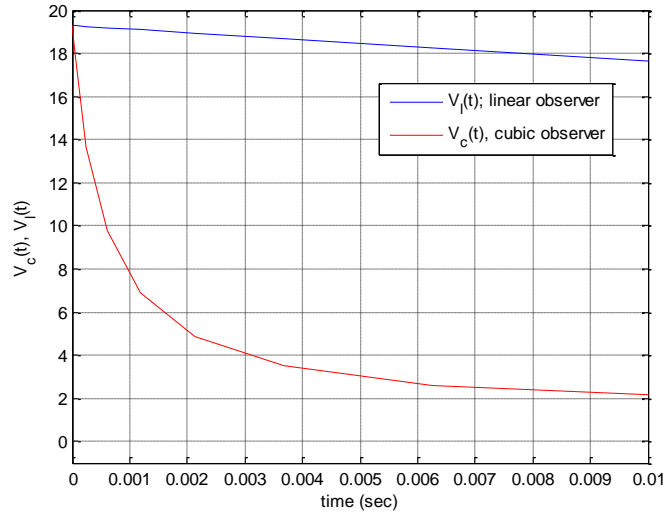

**Figure 2.** Lyapunov functions for estimation errors; the first 10 milliseconds of simulation

According to **Fig. 4** the cubic observer yields a smaller peak (0.48 compared to 1.18) and smaller 5% settling time[٢] (1.7 sec compared to 2.3 sec) for the estimation error. **Fig. 5** compares the estimation error for different values of the scalar $\gamma$. This scalar may be used as a design parameter.

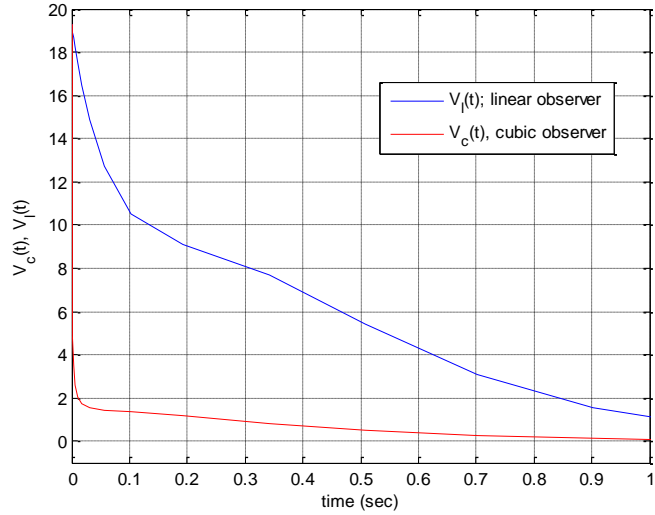

**Figure 3.** Lyapunov functions for estimation error dynamics; the first second of simulation

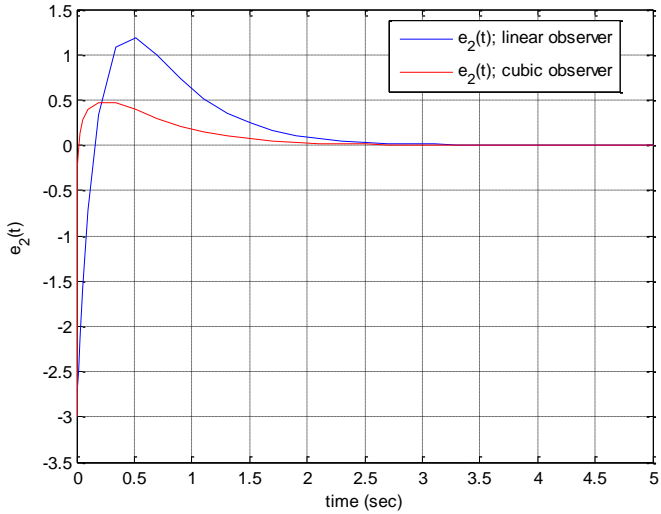

**Figure 4.** Estimation error response for the two observers

---

[٢] We define the settling time as the time during which the error shrinks to less than .05 for the first time.

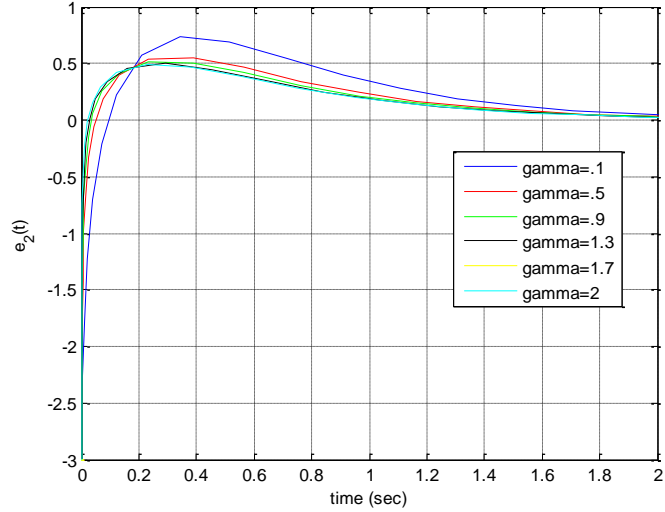

**Figure 5.** Estimation error of the second state for different values of $\gamma$

*Example 2.* This example compares performance and robustness of a cubic observer with a linear one. Consider the following parameter matrices for a linear system:

$$A = \begin{bmatrix} -.1 & -.2 & 0 \\ .3 & 0 & 0 \\ .1 & .2 & -3 \end{bmatrix}, \quad C = [1 \quad 1 \quad 2]$$

Locating observer poles at $[-30 \quad -10 \quad -5]$, the linear observer parameters are derived using pole placement technique. Assuming $Q = 10I_3$, the solution to the Lyapunov equation $P$ is derived as well.

$$L = \begin{bmatrix} 583.7712 \\ -519.9601 \\ -13.4556 \end{bmatrix}, \quad P = \begin{bmatrix} 750.5346 & 785.9444 & 1162.5966 \\ 785.9444 & 823.3524 & 1210.2306 \\ 1162.5966 & 1210.2306 & 2379.6621 \end{bmatrix}, \lambda_{max}(P) = 3702.5756$$

The nonlinear gain $N_c$ is then computed as:

$$N_c = -\gamma P^{-1} C^T \theta = -\gamma \begin{bmatrix} .1866 \\ -.1748 \\ -.0014 \end{bmatrix}$$

**Fig. 6** shows the cumulative squared estimation error for the third state variable for different values of the parameter gamma in (14). The cubic observer with $\gamma = .1$ yields the best performance. The cumulative estimation error norm shown in **Fig. 6** is defined as:

$$J_3(t) = \int_{\tau=0}^{t} (x_3(\tau) - z_3(\tau))^2 d\tau, \quad J(t) = J_1(t) + J_2(t) + J_3(t)$$

In which $z_3(.)$ May be replaced by $\hat{x}_{c3}$ in case of cubic observer and $\hat{x}_{l3}$ in case of linear observer. **Fig. 7** shows the sum of squared estimation errors for all three state variables, $J(t)$ computed for different values of parameter gamma. To evaluate the performances in presence of uncertainties, assume that the state matrix is in the form of (22). (i.e. $A(\varepsilon) = A + \varepsilon I_3$). Recalling *Theorem 5*, the upper bound for $\varepsilon$ is derived as $\varepsilon_{max} \leq .0014$. **Fig. 8** depicts the sum of squared estimation errors for $\varepsilon = .02$ and $\gamma = .1$. The performance advantage persists for $\varepsilon \leq 0.06$. The upper bound .0014, given by *Theorem 5* is very conservative for this case. Tightening of the bounds as well as more generic uncertainty structures require future research.

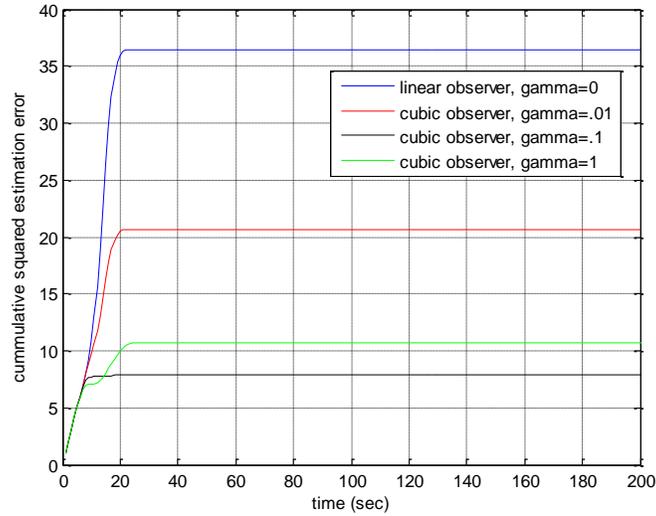

**Figure 6.** Cumulative squared estimation errors of $J_3(t)$, with different values for $\gamma$

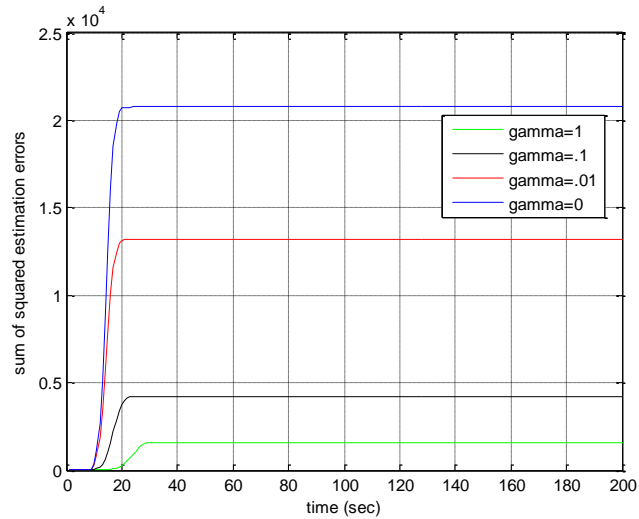

**Figure 7.** Sum of cumulative squared estimation errors $J(t)$, with different values for $\gamma$

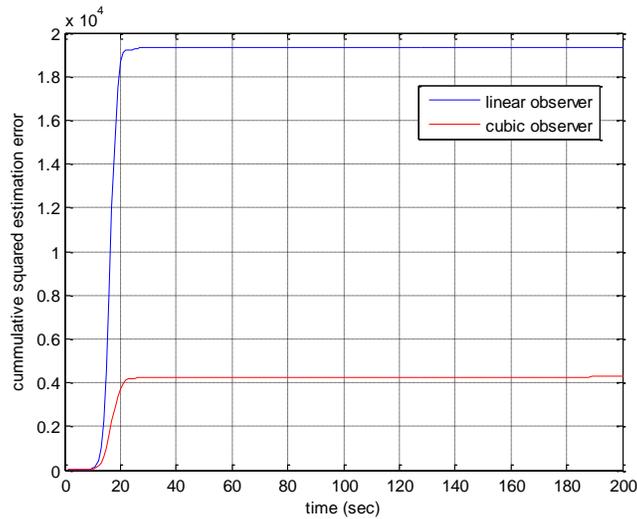

**Figure 8.** Sum of cumulative squared estimation errors $J(t)$, with $\varepsilon = .02, \gamma = .1$

***Example 3.*** This example examines *Theorem 6*. Consider the following system and observer parameters:

$$A = \begin{bmatrix} .1 & -2 & 0 \\ .3 & 0 & -1 \\ .1 & .2 & 3 \end{bmatrix}, \quad C = [1 \ 1 \ 2], \quad B = \begin{bmatrix} 1 & 2 \\ 2 & 0 \\ 0 & 1 \end{bmatrix}, \quad L = L_c = .1 N_c = \begin{bmatrix} .267 \\ -1.429 \\ 3.904 \end{bmatrix},$$

$$N_c = 10L, \quad \theta = 10, \quad K = \begin{bmatrix} -.597 & 2.004 & 2.511 \\ -.197 & .757 & 7.510 \end{bmatrix}$$

**Fig. 9** shows the conventional quadratic cost with $Q_{LQR} = I_3$ and $R_{LQR} = I_2$. Since the feedback gains are the same, the control cost which is the ultimate objective of an observer based control scheme, here, reflects the advantage of the cubic observer over the linear one.

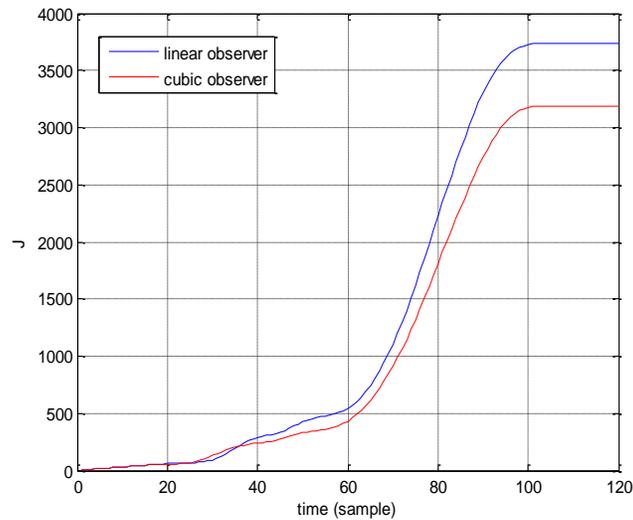

**Figure 9.** Regulation cost by linear observer and cubic observer with linear state feedback

## Conclusion

In this paper, a new nonlinear observer with a cubic term is introduced for improved state estimation of linear systems. It is shown that cubic observers may be used to improve observer performance. These observers can also be readily used for observer-based linear feedback control. The cubic observer is shown to yield better performance in the sense of Lyapunov function decent rate for the estimation error dynamics. It is also shown that the cubic observer has robustness properties similar to linear observers. Examples are included to demonstrate cubic observer performance in comparison to linear observers.